\documentclass{amsart}

\usepackage{amsmath}
\usepackage{amscd}
\usepackage{amssymb}

\input xy
\xyoption{all}

\newcommand{\bC}{{\mathbb C}}

\newcommand{\bN}{{\mathbb N}}
\newcommand{\bP}{{\mathbb P}}

\newcommand{\bR}{{\mathbb R}}
\newcommand{\bZ}{{\mathbb Z}}

\newcommand{\cB}{{\mathcal B}}

\newcommand{\cM}{{\mathcal M}}
\newcommand{\cO}{{\mathcal O}}
\newcommand{\cP}{{\mathcal P}}

\newcommand{\cV}{{\mathcal V}}
\newcommand{\cW}{{\mathcal W}}
\newcommand{\Mbar}{\overline{\cM}}
\newcommand{\vac}{|0\rangle}
\newcommand{\lvac}{\langle 0|}

\DeclareMathOperator{\Aut}{Aut}

\DeclareMathOperator{\val}{val}

\newtheorem{theorem}{Theorem}[section]
\newtheorem{theorem/definition}{Theorem/Definition}[section]

\newtheorem{lemma}{Lemma}[section]

\theoremstyle{remark}

\theoremstyle{definition}

\begin{document}

\title[Localizations and Feynman Rules]
{Localizations on Moduli Spaces and Free Field Realizations of Feynman Rules}
\author{Jian Zhou}
\address{Department of Mathematical Sciences\\Tsinghua University\\Beijing, 100084, China}
\email{jzhou@math.tsinghua.edu.cn} \footnote{This research
is partially supported by research grants from NSFC and Tsinghua
University.}
\begin{abstract}
We prove Iqbal's conjecture on the relationship between the free energy of closed string theory
in local toric geometry and the Wess-Zumino-Witten model.
This is achieved by first reformulating the calculations of the free energy by
localization techniques in terms of suitable Feynman rule,
then exploiting a realization of the Feynman rule by free bosons.
We also use a formula of Hodge integrals conjectured by the author and proved jointly with
Chiu-Chu Melissa Liu and Kefeng Liu.
\end{abstract}
\maketitle

\section{Introduction}

In this work we study an important example of string duality.
Duality in physics literature means the equivalence of different
quantum field theories.
Some examples are already well-known to mathematicians,
e.g. mirror symmetry in string theory.
In recent years,
many more dualities have arisen in string theory.
See e.g. \cite{Vaf} for an exposition.
Most of the string dualities are very mysterious from a mathematical point of view,
they often provide surprising connections among seemingly unrelated mathematical fields.

We will prove a result that mathematically establishes a connection between
the closed string theory of local toric geometry and the Wess-Zumino-Witten theory.
These theories were originally developed by physicists in string theory
for different purposes.
They both have received rigorous mathematical treatments.
The latter corresponds to the representation theory of affine Kac-Moody algebras,
and the former corresponds the intersection theory  of stable moduli spaces,
commonly known as the Gromov-Witten theory.
No connection between these two mathematical theories seems to be previously known
in mathematics literature.

The physicists establish the duality in this case by a sequence of ideas,
including geometric transition,
Chern-Simons theory as string theory,
't Hooft's large $N$ expansion,
relationship between Chern-Simons theory and WZW theory.
Let us give a brief description here.
Suppose a Calabi-Yau three-fold contains a copy of $\bP^1$ with normal bundle $\cO(-1) \oplus \cO(-1)$.
One can performs surgery on  the $\bP^1$,
by replacing it with a copy of $S^3$,
then the resulting space can locally be identified with $T^*S^3$,
the cotangent bundle of $S^3$ (cf e.g. \cite{Gop-Vaf3}).
This process is called the conifold transition in the physics literature.
Witten \cite{Wit} proposed a relationship between the open string theory on the cotangent bundle
of a three-manifold and the Chern-Simons theory on it with structure group $SU(N)$,
by  large $N$ 't Hooft expansion.
Gopakumar and Vafa \cite{Gop-Vaf1} conjectured that under the conifold transition,
the large $N$ Chern-Simons theory on $S^3$ is dual to the $A$-model
closed string theory on $\cO(-1) \oplus(-1) \to \bP^1$.
See also \cite{Gop-Vaf3}.
This conjecture was further tested in \cite{Oog-Vaf}.
By Witten's work \cite{Wit1} on the relationship between Chern-Simons
theory and link invariants,
this duality leads to the idea that both open and string invariants are
related to link invariants.
This has been checked for many cases \cite{Lab-Mar, Ram-Sar, Lab-Mar-Vaf, Mar-Vaf,
Sin-Vaf}.
Closed string invariants for more complicated geometries,
such as local toric del Pezzo surfaces,
have been calculated also from the Chern-Simons theory using geometric transition
\cite{Dia-Flo-Gra1, Dia-Flo-Gra2, Aga-Mar-Vaf}.
Many of these works use
another important idea developed in the seminal paper \cite{Wit1}:
the relationship between Chern-Simons theory and WZW theory.
In mathematics WZW theory is the representation theory of affine Kac-Moody algebras.
See e.g. \cite{Kon2}.
For a fixed integer $k$,
there are only finitely many integrable highest weight representations of level $k$
of an affine Kac-Moody algebra up to equivalence.
Denote their characters by $\chi_0(\tau), \dots, \chi_n(\tau)$.
Then there are holomorphic functions $S_{ij}(\tau)$,
such that
\begin{eqnarray*}
&& \chi_i(- \frac{1}{\tau}) = \sum_j S_{ij}(\tau) \chi_j(\tau).
\end{eqnarray*}
From this one can construction a representation of a double covering of
$SL(2, \bZ)$ (cf. \cite{Ver}).
In this representation,
the  matrix $T$ is a diagonal matrix,
the matrix $S$ is used in the famous Verlinde formula.
A combinatorial description of the matrix $S^{-1}$ for $SU(N)$ in terms of symmetric functions
has been  given in \cite{Mor-Luk, Luk}.

Associate to each very toric Fano surface is a graph obtained
as the image of the moment map of the torus action.
It is a convex polygon where each vertex corresponds to a fixed point,
and each edge corresponds to a one-dimensional orbit of the torus action.
The weight of the induced torus action on the canonical line bundle restricted to a fixed point
determines a ray emanating from the corresponding vertex.
The graph so obtained is called by the physicist the $(p, q)$ five-brane web of the local toric geometry.
In \cite{Aga-Mar-Vaf}
it was shown how to use the $(p, q)$ five-brane web to compute the free energy of closed string theory
in local toric Calabi-Yau geometry via Chern-Simons theory,
and hence the Morton-Lukac formula.
A lattice model interpretation given there inspires Iqbal's conjecture.
He interpreted the $(p, q)$ five-brane web as a Feynman diagram,
and found suitable propagator and vertex to calculate the free energy.
This idea was further developed in \cite{Aga-Kle-Mar-Vaf} where the most general trivalent vertex
are studied and interpreted in terms of open string theory.

It is clear from the above very sketchy description  that it requires a lot of
work to make all the physical arguments involved in the derivation of Iqbal's conjecture
mathematically rigorous.
In this paper,
we present a shortcut to the proof of Iqbal's conjecture,
and in doing so we actually find a unified statement simpler than Iqbal's original conjecture
which is stated case by case (cf. Theorem \ref{thm:Main}).
This is achieved by applying localization techniques on moduli spaces of stable maps.
Our proof relies on the following three key ingredients.
First, we  interpret the terms appearing in localization as Feynman rules.
Secondly, we realize the graphs appearing in localization calculations by free bosonic systems.
This technique was first introduced in \cite{Zho5}.
Thirdly,
we use a formula on Hodge integrals conjectured by the author \cite{Zho4}
and proved jointly with Chiu-Chu Melissa Liu and Kefeng Liu \cite{LLZ}.

Our strategy in this work also sheds some lights on the mathematical understanding of
the topological vertex conjecture \cite{Aga-Kle-Mar-Vaf}.
This will be studied in a forthcoming joint work with Jun Li, Chiu-Chu Melissa Liu and Kefeng Liu.

The rest of the paper is arranged as follows.
In Section \ref{sec:Boson} we recall some definitions and basic facts about free boson systems,
e.g. the vacuum expectation value.
We introduce the notion of abnormal ordering.
In Section \ref{sec:Chemistry} we study a class of labelled graphs called the $\bZ_k$-colored
labelled graphs.
We present a method to create such graphs from a single graph called the $\bZ_k$-cyclic graph
using the free boson systems.
This is referred to as the chemistry of such graphs.
In Section \ref{sec:Feynman}
we use the free boson systems to realize a set of Feynman rules defined for $\bZ_k$-colored
labelled graphs.
We recall in Section \ref{sec:Localization} the interpretation of localization results
in terms of Feynman rule first presented in earlier work of the author \cite{Zho3}.
We prove Iqbal's conjecture in Section \ref{sec:Iqbal}.
Some examples are presented in Section \ref{sec:Examples}.

\section{Free Boson System}
\label{sec:Boson}

We recall in this section some standard results from bosonic string theory.

\subsection{Heisenberg algebra and bosonic Fock space}
The space $\Lambda$ of symmetric functions admits an action of the Heisneberg algebra as follows.
Define operators $\{\beta_n\}_{n \in \bZ}$ on $\Lambda$ by:
\begin{eqnarray*}
&& \beta_{n}(f) = \begin{cases}
p_{-n} f, & n < 0, \\
0, & n =0, \\
n \frac{\partial}{\partial p_{n}}f, & n > 0.
\end{cases}
\end{eqnarray*}
For a partition $\eta$ of length $l$ define:
\begin{align*}
\beta_{\eta} & = \beta_{\eta_1} \cdots \beta_{\eta_l}, &
\beta_{-\eta} & = \beta_{-\eta_1} \cdots \beta_{-\eta_l}.
\end{align*}
Then we have:
\begin{eqnarray*}
&& [\beta_m, \beta_n] = m\delta_{m, -n},\\
&& \beta_n 1 = 0, \;\;\; n \geq 0, \\
&& p_{\eta} = \beta_{-\eta}1.
\end{eqnarray*}
In other words,
$\Lambda$ is the bosonic Fock space in which $1$ is the vacuum vector $\vac$.
Define a Hermitian metric on $\Lambda$ such that
$$\langle p_{\mu}, p_{\nu}\rangle = z_{\mu} \delta_{\mu\nu}.$$
Then in this metric,
one has:
$$\beta_n^* = \beta_{-n},$$
for $n \in \bZ$.

\subsection{Vacuum expectation values and Wick theorem}
Following physicists' notations,
we will write the inner product of $A\vac$ with $\vac$ as
$$\lvac A\vac,$$
where $A$ is a linear operator on $\Lambda$.
It is called the {\em vacuum expectation value} (vev) of $A$,
and will simply be denoted as
$$\langle A \rangle.$$
By the Wick Theorem,
one easily gets:
\begin{eqnarray} \label{eqn:Wick}
&& \langle \beta_{\mu}\beta_{-\nu}\rangle = \delta_{\mu, \nu} z_{\mu},
\end{eqnarray}
where
$$z_{\mu} = \prod_k k^{m_k}m_k!.$$
Here $m_k$ the number of parts of $\mu$ which are equal to $k$.
As a consequence one has:
\begin{eqnarray} \label{eqn:Exponential}
&& \left\langle \exp(\sum_{n \geq 1} \frac{a_nt^n}{n}\beta_n)
\exp(\sum_{n \geq 1} \frac{a_{-n}t^n}{n}\beta_{-n})\right\rangle
= \exp \left(\sum_{n \geq 1} \frac{a_na_{-n}}{n}t^{2n}\right).
\end{eqnarray}

\subsection{Abnormally ordered product}

In physics and mathematics literature,
one usually considers the normally ordered product
$$:\beta_{m_1} \cdots \beta_{m_n}: =
\beta_{m_{i_1}} \cdots \beta_{m_{i_n}},
$$
where $m_{i_1}, \dots, m_{i_n}$ is a permutation of $m_1, \dots, m_n$
such that
$$m_{i_1} \leq \dots \leq m_{i_n}.$$
It is well known that
$$\langle :\beta_{m_1} \cdots \beta_{m_n}:\rangle = 0.$$
We define the {\em abnormally ordered product} as follows:
$$;\beta_{m_1} \cdots \beta_{m_n}; =
\beta_{m_{i_1}} \cdots \beta_{m_{i_n}},
$$
where $m_{i_1}, \dots, m_{i_n}$ is a permutation of $m_1, \dots, m_n$
such that
$$m_{i_1} \geq \dots \geq m_{i_n}.$$
Abnormally ordered products may have nonvanishing vevs,
e.g.
$$\langle ;\beta_m\beta_n;\rangle = |m|\delta_{m, -n}.$$
This will be important to us.

\section{Chemistry of $\bZ_k$-Colored Labelled Graphs}

\label{sec:Chemistry}

\subsection{General definitions and facts about graphs}

For a graph $\Gamma$,
denote by $E(\Gamma)$ the set of edges of $\Gamma$,
$V(\Gamma)$ the set of vertices of $\Gamma$.
The {\em genus} of the graph is given by:
\begin{eqnarray}
&& g(\Gamma) = 1 - |V(\Gamma)| + |E(\Gamma)|. \label{eqn:Genus}
\end{eqnarray}
Recall the valence $val(v)$ of a vertex $v$ is the number of edges incident at $v$.
Since every edge has two vertices,
one has the following identity:
\begin{eqnarray} \label{eqn:Val}
&& \sum_{v \in V(\Gamma)} \val(v) = 2 |E(\Gamma)|.
\end{eqnarray}

\subsection{Cyclic graphs}
\label{sec:Cyclic}

We refer to a graph $\Gamma$ with the following property a {\em cyclic graph}.
There is a one-to-correspondence of between $V(\Gamma)$ and $\bZ_k$,
such that there is exact one edge $e_i$ joining the vertices $v_i$ and $v_{i+1}$ for each $i$,
and there are no other edges.
Denote by $R_n$ the free $\bZ$-module generated by $e_i$, $i\in \bZ_k$,
i.e.,
every element of $R_n$ is a sum:
$$x = \sum_{i \in \bZ_k} a_i e_i$$
for some integers $a_i \in \bZ$.
We associate a system of free bosons to each edge $e_i$:
$$\{\beta_{i, n}: n \in \bZ\}.$$

\subsection{$\bZ_k$-colored labelled graphs}

We refer to a labelled graph with the following property as a {\em $\bZ_k$-colored labelled graph}.
Each vertex $v$ is labelled by an element $i(v) \in \bZ_k$;
each edge $e$ is assigned a natural number $d_e \in \bN$.
Furthermore,
if $v_1$ and $v_2$ are the two vertices of $e$,
then $i(v_1) = i$, $i(v_2) = i+1$, or $i(v_1) = i+1$, $i(v_2) = i$,
and
we write
$$[e] = e_i.$$
We refer to $d_e [e]$ as the degree of the edge.
The {\em degree of the graph} is defined by:
\begin{eqnarray}
&& d(\Gamma) = \sum_{e \in E(\Gamma)} d_e [e] \in R_k. \label{eqn:Degree}
\end{eqnarray}
Denote by $G_g(\bZ_k, d)$ the set of not necessarily connected two-colored labelled graphs of genus $g$
and degree $d$.
The set of connected ones will be denoted by $G_g(\bZ_k, d)^{\circ}$.

\subsection{Chemistry of $\bZ_k$-colored labelled graphs}

For the discussions below,
we introduce some terminologies for $\bZ_k$-colored labelled graphs.
At each vertex with $i(v) =i$,
the edges can be divided into two types depending on the indices of their ends:
if the ends of $e$ have indices $i$ and $i+1$,
then we will say the edge is an outgoing edge;
otherwise, the ends have indices $i$ and $i-1$,
then we will say the edge is an incoming edge.
The degrees of outgoing edges at $v$ determine a partition $\mu^{+}(v)$,
similarly,
the degrees of the incoming edges at $v$ determine a partition $\mu^-(v)$.

Note one of $\mu^+(v)$ and $\mu^-(v)$ might be empty,
but not both.
Denote by $\cP^2_+$ the set of pairs of partitions $(\mu^+, \mu^-)$,
one of which might be empty.
Let $\Gamma \in G_g(\bZ_k, d)$.
For $(\mu^+, \mu^-) \in \cP^2_+$ and $i \in \bZ_k$,
denote by $n^i_{(\mu^+, \mu^-)}(\Gamma)$ the number of vertices $v \in V(\Gamma)$ such that
$i(v) = i$, $\mu^{\pm}(v) = \mu^{\pm}$.

We will refer to $v$ together with the edges incident at it
as an $i(v)$-atom of type $(\mu^+(v), \mu^-(v))$.
The edges can be regarded as chemical bonds that join the atoms.
A $\bZ_k$-colored labelled graph is formed by suitably joining the atoms by the bonds.
As we will show in next section,
this can be realized in bosonic Fock space by Wick theorem.
For that purpose,
we will need the following:

\begin{lemma} \label{lm:GenusDegree}
For any $\Gamma \in G_g(\bZ_k, d)$,
we have:
\begin{eqnarray}
&& \prod_{i \in \bZ_k} \prod_{(\mu^+, \mu^-) \in \cP^2_+} \beta_{i, -\mu^+}^{n^i_{(\mu^+, \mu^-)}(\Gamma)}
= \prod_{i \in \bZ_k} \prod_{(\mu^+, \mu^-) \in \cP^2_+} \beta_{i-1, -\mu^-}^{n^i_{(\mu^+, \mu^-)}(\Gamma)}, \label{eqn:Beta} \\
&& d(\Gamma) = \sum_{i\in\bZ_k} \sum_{(\mu^+, \mu^-) \in \cP^2_+} n^i_{(\mu^+, \mu^-)}(\Gamma) |\mu^+| e_i
= \sum_{i \in \bZ_k} \sum_{(\mu^+, \mu^-) \in \cP_+^2} n^i_{(\mu^+, \mu^-)}(\Gamma) |\mu^-| e_i, \label{eqn:dGamma}  \\
&& 2g(\Gamma) -2 = \sum_{i\in\bZ_k}\sum_{(\mu^+, \mu^-) \in \cP_+^2}
n^i_{(\mu^+, \mu^-)}(\Gamma) (l(\mu^+) + l(\mu^-)- 2).
\end{eqnarray}
\end{lemma}

\begin{proof}
Notice every edge is an outgoing edge for one of its ends,
and an incoming edge for the other end.
This proves  the identity (\ref{eqn:Beta}) and (\ref{eqn:dGamma}).
For each $i$-atom of type $(\mu^+, \mu^-)$,
the valence of the vertex is $l(\mu^+) + l(\mu^-)$.
Hence we have
\begin{eqnarray*}
&& \sum_{i\in \bZ_n} \sum_{(\mu^+, \mu^-) \in \in \cP^2_+} n^i_{(\mu^+, \mu^-)}(\Gamma)
(l(\mu^+) + l(\mu^-) - 2) \\
& = & \sum_{v \in V(\Gamma)} (\val(v) -2)
= \sum_{v \in V(\Gamma)} \val(v) - 2 |V(\Gamma)| \\
& = & 2 |E(\Gamma)| - 2 |V(\Gamma)| = 2 g(\Gamma) - 2.
\end{eqnarray*}
Here in the last equality we have used (\ref{eqn:Genus}).
\end{proof}

\subsection{Automorphism groups of $\bZ_k$-colored labelled graphs}

An automorphism of a $\bZ_k$-colored labelled graph $\Gamma$
consists of two one-to-one correspondences:
$f^V: V(\Gamma) \to V(\Gamma)$ and $f^E: E(\Gamma) \to E(\Gamma)$,
with the following requirements.
\begin{itemize}
\item $i(f^V(v)) = i(v)$, for all $v \in V(\Gamma)$;
\item $d_{f^E(e)} = d_e$, for all $e \in E(\Gamma)$;
\item if $v_1, v_2 \in V(\Gamma)$ are the two vertices of an edge $e \in E(\Gamma)$,
then $f^V(v_1)$ and $f^V(v_2)$ are the two vertices of $f^E(e)$.
\end{itemize}

\begin{lemma} \label{lm:Automorphism}
Suppose $\{n^i_{(\mu^+, \mu^-)}: i\in \bZ_k, (\mu^+, \mu^-) \in \cP^2_+\}$
is a collection of nonnegative integers,
which contain only finitely many nonzero integers and satisfy (\ref{eqn:Beta}).
Then we have:
\begin{eqnarray*}
&& \left\langle ;\prod_{i\in \bZ_k} \prod_{(\mu^+, \mu^-) \in \cP^2}
\frac{1}{n^i_{(\mu^+, \mu^-)}!} \left(\frac{\beta_{i, \mu^+}}{z_{\mu^+}}\cdot
\frac{\beta_{i-1, -\mu^-}}{z_{\mu^-}}\right)^{n^i_{(\mu^+, \mu^-)}};\right\rangle\\
& = & \sum_{\Gamma} \frac{1}{|\Aut_{\Gamma}| \cdot \prod_{e \in E(\Gamma)} d_e},
\end{eqnarray*}
where the sum is taken over all graphs $\Gamma \in G_g(\bZ_k, d)$ which
satisfy:
\begin{eqnarray*}
&& n^i_{(\mu^+, \mu^-)}(\Gamma) = n^i_{(\mu^+, \mu^-)}.
\end{eqnarray*}
In the exceptional case when all $n^i_{(\mu^+, \mu^-)} = 0$,
the vev is $1$.
\end{lemma}

\begin{proof}
This is an easy consequence of the Wick Theorem.
\end{proof}

\section{Generalized Vertex Operators and Feynman Rule}
\label{sec:Feynman}

\subsection{Generalized vertex operators}
Suppose we are given a collection
$$\{w_{i, (\mu^+_i,\mu_i^-)}: i \in \bZ_k, (\mu^+_i,\mu_i^-) \in \cP_+^2\}.$$
For each $i \in \bZ_k$,
define
\begin{eqnarray*}
&& Y_i(\beta) = \sum_{(\mu^+_i, \mu_i^-) \in \cP^2_+}
w_{i, (\mu^+_i,\mu_i^-)}\frac{\beta_{i, \mu^+_i}}{z_{\mu^+_i}} \lambda^{l(\mu^+_i)-1}
t_i^{\frac{|\mu^+_i|}{2}}
\cdot \frac{\beta_{i-1, -\mu^-_i}}{z_{\mu_i^-}}\lambda^{l(\mu^-_i)-1}
t_{i-1}^{\frac{|\mu^-_i|}{2}}, \\
&& X_{i}(\beta) = e^{Y_i(\beta)}.
\end{eqnarray*}
In the following,
$t_i$ will be referred to as the {\em degree tracking parameter},
$\lambda$ will be referred to as the {\em genus tracking parameter}.
We will consider the correlation function
\begin{eqnarray*}
&& Z = \langle ; \prod_{i \in \bZ_k} X_i(\beta);\rangle,
\end{eqnarray*}
and the free energy:
$$F = \log Z.$$

\subsection{Free field realizations of Feynman rules}

Now we regard $X_i(\beta)$ and $Y_i(\beta)$ as collections of $\pm$-atoms.
Taking the vevs
\begin{eqnarray*}
&& \langle ;\prod_{i \in \bZ_k} X_i(\beta);\rangle
\end{eqnarray*}
can be regarded as considering all chemical reactions in which the atoms are joined
by the chemical bonds.
This is exactly the context of Wick theorem.
More precisely,
we have the following result.
Before stating it,
let us introduce some notations.
For $d = (d_1, \dots, d_k)$,
we write $d \geq 0$ if $d_1, \dots, d_k \geq 0$;
$d > 0$ if $d \geq 0$ and one of $d-1, \dots, d_k$ is positive.
We also write
$$t^d = t_1^{d_1} \cdots t_k^{d_k}.$$

\begin{theorem} \label{thm:Feynman}
The following identities hold:
\begin{eqnarray}
Z & = & \sum_{d \geq 0} \sum_{g \geq 0} \lambda^{2g-2} \sum_{\Gamma \in G_g(\bZ_k, d)}
\frac{1}{|A_{\Gamma}|} \prod_{v \in V(\Gamma)} w_v \prod_{e \in E(\Gamma)} w_e,
\label{eqn:Z} \\
F & = & \sum_{d > 0} \sum_{g \geq 0} \lambda^{2g-2} \sum_{\Gamma \in G_g(\bZ_k, d)^{\circ}}
\frac{1}{|A_{\Gamma}|} \prod_{v \in V(\Gamma)} w_v \prod_{e \in E(\Gamma)} w_e. \label{eqn:F}
\end{eqnarray}
Here for a vertex $v$ of type $(\eta^+, \eta^-)$ with $i(v) = i$,
\begin{align*}
w_v = w_{i, (\eta^+, \eta^-)};
\end{align*}
for an edge $e$ with $i(e) = i$,
$$w_e = t_{i}^{d_e}.$$
\end{theorem}

\begin{proof}
The identity (\ref{eqn:Z}) is an easy consequence of the Wick theorem.
\begin{eqnarray*}
&& Z  = \langle ;\prod_{i\in\bZ_k} X_i(\beta); \rangle \\
& = & \left\langle ;\prod_{i \in \bZ_k}\prod_{(\mu^+, \mu^-) \in \cP_+^2}
\exp\left(w_{i, (\mu^+, \mu^-)}\frac{\beta_{i, \mu^+}}{z_{\mu^+}}
\lambda^{l(\mu^+)-1} t_i^{\frac{|\mu^+|}{2}}
\frac{\beta_{i-1, -\mu^-}}{z_{\mu^-}} \lambda^{l(\mu^-)-1}
t_{i-1}^{\frac{|\mu^-|}{2}}\right);\right\rangle \\
& =&  \left\langle ;\prod_{i \in \bZ_k} \prod_{(\mu^+, \mu^-) \in \cP^2_+}
\sum_{n^i_{\mu^+, \mu^-}} \frac{\left(w_{i, (\mu^+, \mu^-)}
\frac{\beta_{i, \mu^+}}{z_{\mu^+}}
\cdot \frac{\beta_{i-1, -\mu^-}}{z_{\mu^-}}
\lambda^{l(\mu^+)+l(\mu^-)-2}t_i^{\frac{|\mu^+|}{2}}
t_{i-1}^{\frac{|\mu^-|}{2}}\right)
^{n^i_{(\mu^+, \mu^-)}}}{n^i_{\mu^+, \mu^-}!}; \right\rangle\\
& = & \sum_{\{n^i_{(\mu^+, \mu^-)}\}}
\lambda^{\sum_{i \in \bZ_k}\sum_{(\mu^+, \mu^-)\in \cP_+^2}
\left[(l(\mu^+) + l(\mu^-) - 2)n^i_{(\mu^+, \mu^-)}\right]} \\
&& \cdot \prod_{i \in \bZ_k}
t_i^{\frac{1}{2} \sum_{(\mu^+, \mu^-)\in \cP_+^2} |\mu^+|n^i_{(\mu^+, \mu^-)}}
t_{i-1}^{\frac{1}{2} \sum_{(\mu^+, \mu^-)\in \cP_+^2} |\mu^-|n^i_{(\mu^+, \mu^-)}}\\
&& \cdot \prod_{i \in \bZ_k} \prod_{(\mu^+, \mu^-) \in \cP_+^2}
w_{i, (\mu^+, \mu^-)}^{n^i_{(\mu^+, \mu^-)}} \cdot
\left\langle ;\prod_{i \in \bZ_k} \prod_{(\mu_+, \mu_-) \in \cP_+^2}
\frac{;\left(\frac{\beta_{i, \mu^+}}{z_{\mu^+}}\cdot
\frac{\beta_{i-1, -\mu^-}}{z_{\mu^-}} \right)^{n^i_{(\mu^+, \mu^-)}}}{n^i_{(\mu^+, \mu^-)}!} ;\right\rangle\\
& = & \sum_{d} \sum_{g \geq 0} \sum_{\Gamma \in G_g(\bZ_k, d)}
\frac{1}{|A_{\Gamma}|} \lambda^{2g(\Gamma)-2} t^d
\prod_{i \in \bZ_k} \prod_{\mu}  w_{i, (\mu^+, \mu^-)}^{n^i_{(\mu^+, \mu^-)}(\Gamma)}.
\end{eqnarray*}
In the last equality we have used Lemma \ref{lm:GenusDegree} and Lemma \ref{lm:Automorphism}.
This proves (\ref{eqn:Z}).

It is a standard result that the exponential of the right side of (\ref{eqn:F})
is the right-hand side of (\ref{eqn:Z}).
This proves (\ref{eqn:F}).
\end{proof}

In the following,
we will obtain Feynman rules by localizations on moduli spaces of stable maps,
then apply Theorem \ref{thm:Feynman} to give free field realization of the Feynman rules.
Then standard techniques for bosonic system can be applied to carry out the calculations.
%We will give the details in the cases of local $\bP^1 \times \bP^1$ geometries.

\section{Localization on Moduli Spaces of Stable Maps to Surfaces}
\label{sec:Localization}

In this section we recall the results in \cite{Zho3}
on the localizations on the moduli space $\Mbar_{g, 0}(X, d)$
for a toric Fano surface $X$.

\subsection{Toric Fano surfaces and cyclic graphs}

Let $X$ be a toric Fano surface with associated $T = (\bC^*)^2$-action.
The image $\mu(X)$ of the moment map $\mu: X \to \bR^2$ of the $T$-action is a convex polygon
whose vertices are the images of the fixed points  (see e.g. \cite{Ati}).
At each fixed point $p$ of $X$,
there is a decomposition
$$T_pX = L^1_{p} \oplus L^2_{p},$$
where $L^j_{p}$ are one-dimensional subspaces on which $T$ acts with weight $\lambda_{pj}$.
These two weights are linearly independent,
each corresponding to a one-dimensional orbit of the $T$-action.
The closure of each of these orbits is a copy of $\bP^1$,
giving an equivariant map $f_{pq}: \bP^1 \to X$,
such that
\begin{align*}
f_{pq}([0:1]) & = p, & f_{pq}([1:0]) & = q,
\end{align*}
where $q$ is another fixed point of $X$.
The image of $f_{pq}(\bP^1)$ under the moment map is exactly the edge of $\mu(X)$ joining $\mu(p)$ to $\mu(q)$.
Note both $X$ and $T$ have canonical orientations,
hence $\mu(S)$ has an induced orientation.
We regard the polygon $\mu(X)$ as a cyclic graph and give its vertices and edges
the labelling as in \S \ref{sec:Cyclic} in accordance with the induced orientation.
Denote the weights of $T_{p_{i(v)}}X$ by $u_{i(v), i(v)+1}$ and $u_{i(v), i(v)-1}$.
Then clearly:
$$u_{i(v), i(v)+1} = -u_{i(v)+1, i(v)}.$$

\subsection{Fixed points on $\Mbar_{g, 0}(X, d)$}

For $d \in H_2(X, \bZ)$,
denote by $\Mbar_{g, 0}(X, d)$ the moduli space of stable maps of genus $g$
to $X$ of class $d$.
The $T$-action induces $T$-actions on $\Mbar_{g,0}(X, d)$.
The fixed point components of $\Mbar_{g, n}(X, d)^T$ are very easy to describe.
They are in one-to-one correspondence with a set $G_g(X, d)$ of
decorated graphs described below.
Each vertex $v$ of the graph $\Gamma \in G_g(X, d)$ is assigned an index $i(v) \in X^T$,
and a genus $g(v)$.
The valence $\val(v)$ of $v$ is the number of edges incident at $v$.
If two vertices $u$ and $v$ are joined by an edge $e$,
then $i(u) \neq i(v)$,
and $e$ is assigned a ``degree"
$$\delta(e) = d_e [l_{i(u)i(v)}] \in H_2(X, \bZ).$$
Denote by $E(\Gamma)$ the set of edges of $\Gamma$,
$V(\Gamma)$ the set of vertices of $\Gamma$.
The genus of the graph is given by
$$g(\Gamma) = 1 - |V(\Gamma)| +|E(\Gamma)|.$$
The decorations of $\Gamma$ are required to satisfy the following conditions:
\begin{align*}
\sum_{e \in E(\Gamma)} \delta_e & = d, &
\sum_{v \in V(\Gamma)} g(v) + g(\Gamma) & = g.
\end{align*}
Let $f: C \to X$ represent a fixed point.
Then each vertex $v$ corresponds to a connected component $C_v$ of genus $g(v)$,
with $\val(v)$ nodal points.
The component $C_v$ is mapped by $f$ to the fixed point $i(v)$.
When $2g(v) -2  + \val(v) < 0$,
$C_v$ is simply a point.
There are only two cases when this happens:
$g(v) = 0$ and $\val(v) = 1$, $g(v) = 0$ and $\val(v) = 2$.
They will be referred to as the type I and type II unstable vertices respectively.
Each edge $e$ corresponds to a component of $C$,
isomorphic to $\bP^1$.
Each $C_e$ is mapped to the balloon $l_{i(v)i(u)}$ with degree $d_e$.

A flag $F$ is a pair $(v, e)$,
where $v$ is a vertex and $e$ is an edge incident at $v$.
Each flag $F = (v, e)$ corresponds to a nodal point $x_F$ of $C$.
Denote by $L_F$ the line bundle on $\Mbar_{\Gamma}$ whose fiber at
$f: C \to X$ is give by the cotangent space to $C_v$ at $x_F$.
Denote by $\psi_F$ the first Chern class of $L_F$.

Define
$$\Mbar_{\Gamma} = \prod_{v \in V(\Gamma)} \Mbar_{g(v), val(v)}.$$
In this product,
$\Mbar_{0,1}$ and $\Mbar_{0, 2}$ are interpreted as points.
There are natural morphisms
$$\tau_{\Gamma}: \Mbar_{\Gamma} \to \Mbar_{g,0}(\bP^1, d).$$
Its image is $\Mbar_{\Gamma}/A_{\Gamma}$,
where for $A_{\Gamma}$ we have an exact sequence:
$$0 \to \prod_{e \in E(\Gamma)} \bZ_{d_e} \to A_{\Gamma} \to \Aut(\Gamma) \to 1.$$

Given a graph $\Gamma$ in $G_g(X, d)$,
we call the labelled graph obtained from $\Gamma$ by ignoring
the markings of $g(v)$ of the vertices the type of $\Gamma$.
Denote by $G(X, d)$ the set of types of graphs in $G_{g}(X, d)$.

\subsection{Localization on moduli spaces}

Introduce the following notation.
Suppose $(v, e)$ is a flag,
then
$$u_{(v, e)} = u_{i(v), i(v) \pm 1},$$
where the $\pm$ sign depends on whether $e$ is an outgoing or incoming edge at $v$.

Note $TX|_{f(C_e)}$ decompose into the a direct sum of the tangent bundle and the normal bundle of $f(C_e)$.
The tangent bundle has weights $u_{i(v), i(v)+1}$ at $p_{i(v)}$,
$u_{i(v)+1, i(v)}$ at $p_{i(v)+1}$;
the normal bundle has weights $u_{i(v), i(v)-1}$ at $p_{i(v)}$,
$u_{i(v)+1, i(v)+2}$ at $p_{i(v)+1}$.
Now the tangent bundle has degree
$$\frac{u_{i(v), i(v)+1} -  u_{i(v)+1, i(v)}}{u_{i(v),i(v)+1}} = 2,$$
while the normal bundle has degree
\begin{eqnarray} \label{eqn:s}
&& s_e:=\frac{u_{i(v), i(v)-1} -  u_{i(v)+1, i(v)+2}}{u_{i(v), i(v)+1}}.
\end{eqnarray}
By adjunction formula,
$$2 + s_e = - K_X \cdot f(C_e) > 0,$$
hence
$$s_e > -2.$$
Finally, note the weight of $K_X$ at $i(v)$ is
$$- u_{i(v), i(v)+1} - u_{i(v), i(v)-1}.$$

Now by \cite{Zho3} we have
\begin{eqnarray} \label{eqn:Localization}
\int_{[\Mbar_{g,0}(X, d)]^{vir}} e_T(U_d^g)
& = & \sum_{\Gamma \in G_{g}(X, d)} \frac{1}{|A_{\Gamma}|}
\int_{\Mbar_{\Gamma}} \frac{\tau^*_{\Gamma}(e_T(U_d^g))}{e_T(N_{\Gamma})},
\end{eqnarray}
where
\begin{eqnarray}
&& \int_{\Mbar_{\Gamma}} \frac{\tau^*_{\Gamma}(e_T(U_d^g))}{e_T(N_{\Gamma})}
= \prod_{v \in V(\Gamma)} w_v \cdot \prod_{e \in E(\Gamma)} w_e.
\end{eqnarray}
For $v \in V^s(\Gamma)$,
\begin{equation} \label{eqn:wv}
\begin{split}
w_v = & \int_{\Mbar_{g(v)}, \val(v)}
\frac{\Lambda_{g(v)}^{\vee}(u_{i(v),i(v)+1})\Lambda^{\vee}_{g(v)}(u_{i(v), i(v)-1})
\Lambda_{g(v)}^{\vee}(-u_{i(v), i(v)+1} - u_{i(v),i(v)-1})}
{\prod_{(v, e) \in F(\Gamma)} \left(\frac{u_{(v, e)}}{d_e} - \psi_{(v, e)}\right)} \\
& \cdot \left[u_{i(v), i(v)+1}u_{i(v), i(v)-1}(-u_{i(v), i(v)+1} - u_{i(v), i(v)-1})\right]^{\val(v)-1};
\end{split} \end{equation}
for $v \in V^I(\Gamma)$,
\begin{eqnarray}
w_v & = & \frac{u_{(v, e)}}{d_e},
\end{eqnarray}
where $(v, e) \in F^I(\Gamma)$;
for $v \in V^{II}(\Gamma)$,
\begin{eqnarray}
w_v & = & \frac{u_{i(v), i(v)+1}u_{i(v), i(v)-1}(- u_{i(v), i(v)+1} - u_{i(v), i(v)-1})}
{\frac{u_{(v, e_1)}}{d_{e_1}} + \frac{u_{(v, e_2)}}{d_{e_2}}},
\end{eqnarray}
where $e_1$ and $e_2$ are the two edges incident at $v$.
For $e \in E(\Gamma)$,
\begin{eqnarray*}
w_e & = & -1 \cdot (-1)^{s_ed_e} \cdot  \frac{\prod_{a=1}^{d_e -1}
\left(d_e u_{i(v), i(v)-1} + a u_{i(v), i(v)+1}\right)}
{\frac{d_e!}{d_e}u_{i(v), i(v)+1}^{d_e}} \\
&& \cdot  \frac{\prod_{a=1}^{d_e -1}
\left(d_e u_{i(v)+1, i(v)+2} + a u_{i(v)+1, i(v)}\right)}
{\frac{d_e!}{d_e}u_{i(v)+1, i(v)}^{d_e}}.
\end{eqnarray*}
We will reformulate these rules in the next section.

\section{Proof of Iqbal's Conjecture}
\label{sec:Iqbal}

\subsection{A formula for Hodge integrals}

Consider the following generating series of Hodge integrals:
\begin{eqnarray*}
&& G_{\mu^+, \mu^-}(x, y) \\
& = & - \frac{\sqrt{-1}^{l(\mu^+)+l(\mu^-)}}{z_{\mu^+} \cdot z_{\mu^-}}
\sum_{g \geq 0} \lambda^{2g-2+l(\mu^+)+l(\mu^-)} \\
&& \int_{\Mbar_{g, l(\mu^+)+l(\mu^-)}}
\frac{\Lambda_{g}^{\vee}(x)\Lambda^{\vee}_{g}(y)\Lambda_{g}^{\vee}(-x - y)}
{\prod_{i=1}^{l(\mu^+)} \frac{x}{\mu_i^+} \left(\frac{x}{\mu^+_i} - \psi_i\right)
\prod_{j=1}^{l(\mu^-)} \frac{y}{\mu_i^-}\left(\frac{y}{\mu^-_j} - \psi_{j+l(\mu^+)}\right)} \\
&& \cdot \left[xy(x+y)\right]^{l(\mu^+)+l(\mu^-)-1}
\cdot \prod_{i=1}^{l(\mu^+)} \frac{\prod_{a=1}^{\mu^+_i-1}
\left( \mu^+_iy + a x\right)}{\mu_i^+! x^{\mu_i^+-1}}
\cdot \prod_{i=1}^{l(\mu^-)} \frac{\prod_{a=1}^{\mu^-_i-1}
\left( \mu_i^- x + a y\right)}{\mu_i^-! y^{\mu_i^--1}}.
\end{eqnarray*}
The foloowing formula was conjectured by the author in \cite{Zho4}.
See \cite{LLZ} for its proof.
\begin{equation} \label{eqn:Conjecture}
\begin{split}
& \exp \left(
\sum_{(\mu^+, \mu^-) \in \cP_+^2} G_{\mu^+, \mu^-}(\tau)p^+_{\mu^+}p^-_{\mu^-}\right) \\
= & \sum_{|\mu^{\pm}|=|\nu^{\pm}|}
\frac{\chi_{\nu^+}(\mu^+)}{z_{\mu^+}} \frac{\chi_{\nu^-}(\mu^-)}{z_{\mu^-}}
e^{\sqrt{-1}(\kappa_{\nu^+} \tau  + \frac{\kappa_{\nu^-}}{\tau})\lambda/2}
\cW_{\nu^+, \nu^-} p^+_{\mu^+}p^-_{\mu^-}.
\end{split} \end{equation}

\subsection{Reformulation of the Feynman rule for a graph}

Now we reformulate the Feynman rule for a graph in $G_{g,0}(X, d)$.
First of all,
we rewrite the vertex as
\begin{eqnarray*}
&& w_v \\
& =&- \sqrt{-1}^{l(\mu^+(v))+l(\mu^-(v))}G_{\mu^+(v), \mu^-(v)}(u_{i(v)}^+, u_{i(v)}^{-1}) \\
&& \cdot [u_{i(v)}^+u_{i(v)}^-(u_{i(v)}^++u_{i(v)}^-)]^{l(\mu^+(v))+l(\mu^-(v))-1} \\
&& \cdot \sqrt{-1}^{l(\mu^+(v))+l(\mu^-(v))} \prod_{j=1}^{l(\mu^+(v))} \frac{u_{i(v)}^+}{\mu^+_j(v)}
\cdot \prod_{j=1}^{l(\mu^-(v))} \frac{u_{i(v)}^-}{\mu^-_j(v)}. \nonumber
\end{eqnarray*}
We will deal with the four factors in $w_e$ as follows.
The third and the fourth factor will be reassigned to the corresponding vertices;
the second factor $(-1)^{s_ed_e}$ remains;
the first factor $-1$ cancels with a factor in $w_v$ as follows:
\begin{eqnarray*}
&& \prod_{v \in \Gamma} \sqrt{-1}^{l(\mu^+(v))+l(\mu^-(v)}
\cdot \prod_{e \in E(\Gamma)} (-1) \\
& = & \prod_{v \in V(\Gamma)} \sqrt{-1}^{\val(v)} \cdot (-1)^{|E(\Gamma)|} \\
& = & \sqrt{-1}^{\sum_{v\in V(\Gamma)} \val(v)} \cdot (-1)^{|E(\Gamma)|} \\
& = & \sqrt{-1}^{2|E(\Gamma)|} \cdot (-1)^{|E(\Gamma)|} = 1.
\end{eqnarray*}
Hence we can get a new Feynman rule as follows:
\begin{eqnarray*}
w_v & = & - \sqrt{-1}^{l(\mu^+(v))+l(\mu^-(v))}G_{\mu^+(v), \mu^-(v)}(u_{i(v)}^+, u_{i(v)}^{-1}) \\
&& \cdot [u_{i(v)}^+u_{i(v)}^-(u_{i(v)}^++u_{i(v)}^-)]^{l(\mu^+(v))+l(\mu^-(v))-1} \\
&& \cdot \prod_{i=1}^{l(\mu^+(v))} \frac{\prod_{a=1}^{\mu^+_i(v)-1}
\left( \mu^+_i(v)y + a x\right)}{\mu_i^+(v)! x^{\mu_i^+(v)-1}}
\cdot \prod_{i=1}^{l(\mu^-(v))} \frac{\prod_{a=1}^{\mu^-_i(v)-1}
\left( \mu_i^-(v) x + a y\right)}{\mu_i^-(v)! y^{\mu_i^-(v)-1}},
\end{eqnarray*}
and
\begin{eqnarray*}
&& w_e = (-1)^{s_ed_e}.
\end{eqnarray*}

\subsection{Feynman rule for a type of graphs}

Now we consider
\begin{eqnarray*}
F_d(\lambda) & = & \sum_{g \geq 0} \lambda^{2g-2} \int_{[\Mbar_{g,0}(X, d)]^{vir}} e_T(U_d^g) \\
& = &  \sum_{g \geq 0} \lambda^{2g-2}  \sum_{\Gamma \in G_{g}(X, d)} \frac{1}{|A_{\Gamma}|}
\int_{\Mbar_{\Gamma}} \frac{\tau^*_{\Gamma}(e_T(U_d^g))}{e_T(N_{\Gamma})},
\end{eqnarray*}
and its generating series:
$$F(\lambda) = \sum_{d \in H_2(X, \bZ) - \{0\}} F_d(\lambda) t^d.$$
By the result of preceding subsection,
it is straightforward to get
\begin{eqnarray*}
F(\lambda) & = & \sum_{g \geq 0} \lambda^{2g-2} \sum_{\Gamma \in G_g(\bZ_k, d)}
\frac{1}{|A_{\Gamma}|} \prod_{v \in V(\Gamma)} w_v \cdot \prod_{e \in E(\Gamma)} w_e,
\end{eqnarray*}
where
\begin{eqnarray*}
&& w_v = z_{\mu^+(v)} \cdot z_{\mu^-(v)} \cdot
G_{\mu^+(v), \mu^-(v)}(u_{i(v)}^+, u_{i(v)}^-), \\
&& w_e = ((-1)^{s_e}t_{i(e)})^{d_e}.
\end{eqnarray*}
Here
$$t_i = e^{[l_{i, i+1}]}.$$

From now on we will use the following notations:
\begin{align*}
u_i^+ & = u_{i, i+1}, & u_i^- & = u_{i, i-1}.
\end{align*}
Note we have
\begin{eqnarray} \label{eqn:si}
&& \frac{u^-_i}{u^+_i} + \frac{u_{i+1}^+}{u_{i+1}^-}
= \frac{u^-_i - u^+_{i+1}}{u_i^+} = s_i.
\end{eqnarray}

\subsection{Proof of Iqbal's conjecture}
Here we present a unified statement of Iqbal's conjecture which was originally stated case by case
in \cite{Iqb}.

\begin{theorem} \label{thm:Main}
Assuming (\ref{eqn:Conjecture}),
we have
\begin{eqnarray}
&& Z(\lambda)
= \prod_{i\in \bZ_k} \sum_{\nu_i}
e^{\sqrt{-1} \kappa_{\nu_i}s_i \lambda/2}
\cW_{\nu_i, \nu_{i-1}} ((-1)^{s_i}t_i)^{|\nu_i|}.
\end{eqnarray}
\end{theorem}

\begin{proof}
Combining the Feynman rule above with Theorem \ref{thm:Feynman},
we have
\begin{eqnarray*}
Z(\lambda) & = & \langle ; \prod_{i\in\bZ_k} e^{Y_i(\beta)};\rangle,
\end{eqnarray*}
where
\begin{eqnarray*}
&& Y_i(\beta) \\
& = & \sum_{(\mu^+_i, \mu_i^-) \in \cP^2_+}
w_{i, (\mu^+_i,\mu_i^-)}\frac{\beta_{i, \mu^+_i}}{z_{\mu^+_i}} \lambda^{l(\mu^+_i)-1}
(\sqrt{-1}^{s_i}t_i)^{\frac{|\mu^+_i|}{2}} \\
&& \cdot \frac{\beta_{i-1, -\mu^-_i}}{z_{\mu_i^-}}\lambda^{l(\mu^-_i)-1}
(\sqrt{-1}^{s_{i-1}}t_{i-1})^{\frac{|\mu^-_i|}{2}},
\end{eqnarray*}
for
$$w_{i, (\mu_i^+, \mu^-_i)}
=  z_{\mu^+_i} \cdot z_{\mu^-_i} \cdot
G_{\mu^+_i, \mu^-_i}(u_{i}^+, u_{i}^-).$$
Hence
\begin{eqnarray*}
Y_i(\beta) & = &
 \sum_{(\mu^+_i, \mu_i^-) \in \cP^2_+}
G_{\mu^+_i,\mu_i^-}(u_i^+, u_i^-) \beta_{i, \mu^+_i} \lambda^{l(\mu^+_i)-1}
((-1)^{s_i}t_i)^{\frac{|\mu^+_i|}{2}} \\
&& \cdot \beta_{i-1, -\mu^-_i}\lambda^{l(\mu^-_i)-1}
((-1)^{s_{i-1}}t_{i-1})^{\frac{|\mu^-_i|}{2}},
\end{eqnarray*}
and so by (\ref{eqn:Conjecture}),
\begin{eqnarray*}
&& X_i(\beta) = e^{Y_i(\beta)} \\
& = & \sum_{|\mu_i^{\pm}|=|\nu_i^{\pm}|}
\frac{\chi_{\nu_i^+}(\mu_i^+)}{z_{\mu_i^+}} \frac{\chi_{\nu_i^-}(\mu_i^-)}{z_{\mu_i^-}}
e^{\sqrt{-1}\left(\kappa_{\nu^+} \frac{u^-_i}{u^+_i} + \kappa_{\nu^-}\frac{u_i^+}{u_i^-}\right)\lambda/2}
\cW_{\nu_i^+, \nu_i^-} \\
&& \cdot \beta_{i, \mu_i^+} ((-1)^{s_i}t_i)^{\frac{|\mu^+_i|}{2}}
\cdot \beta_{i-1, -\mu_i^-}((-1)^{s_{i-1}}t_{i-1})^{\frac{|\mu^-_i|}{2}}.
\end{eqnarray*}
It follows that
\begin{eqnarray*}
Z & = & \langle ; \prod_{i \in \bZ_k} X_i(\beta) ;\rangle \\
& = & \left\langle ;\prod_{i \in \bZ_k}\sum_{|\mu_i^{\pm}|=|\nu_i^{\pm}|}
\frac{\chi_{\nu_i^+}(\mu_i^+)}{z_{\mu_i^+}} \frac{\chi_{\nu_i^-}(\mu_i^-)}{z_{\mu_i^-}}
e^{\sqrt{-1}\left(\kappa_{\nu^+} \frac{u^-_i}{u^+_i} + \kappa_{\nu^-}\frac{u_i^+}{u_i^-}\right)\lambda/2}
\cW_{\nu_i^+, \nu_i^-} \right.\\
&& \left.\cdot \beta_{i, \mu_i^+} ((-1)^{s_i}t_i)^{\frac{|\mu^+_i|}{2}}
\cdot \beta_{i-1, -\mu_i^-}((-1)^{s_{i-1}}t_{i-1})^{\frac{|\mu^-_i|}{2}};\right\rangle.
\end{eqnarray*}
Note the nonzero contributions come from terms with
$$\mu_i^+ = \mu_{i+1}^-,$$
$i \in \bZ_k$.
Hence
\begin{eqnarray*}
&& Z(\lambda) \\
& = & \prod_{i\in \bZ_k} \sum_{\mu_i^+, \nu_i^{\pm}}
\frac{\chi_{\nu_i^+}(\mu_i^+)}{z_{\mu_i^+}}
\frac{\chi_{\nu_i^-}(\mu_{i-1}^+)}{z_{\mu_{i-1}^+}}
e^{\sqrt{-1}\left(\kappa_{\nu_i^+} \frac{u^-_i}{u^+_i} + \kappa_{\nu_i^-}\frac{u_i^+}{u_i^-}\right)\lambda/2}
\cW_{\nu_i^+, \nu_i^-} \\
&& \cdot z_{\mu_i^+} ((-1)^{s_i}t_i)^{|\mu^+_i|} \\
& = & \prod_{i\in \bZ_k} \sum_{\mu_i^{+}, \nu_i^{\pm}}
\frac{\chi_{\nu_i^+}(\mu_i^+)\chi_{\nu_i^-}(\mu_{i-1}^+)}{z_{\mu_{i-1}^+}}
e^{\sqrt{-1}\left(\kappa_{\nu_i^+} \frac{u^-_i}{u^+_i} + \kappa_{\nu_i^+}\frac{u_i^+}{u_i^-}\right)\lambda/2}
\cW_{\nu_i^+, \nu_i^-} ((-1)^{s_i}t_i)^{|\nu^+_i|}.
\end{eqnarray*}
Now for each $i \in \bZ_k$,
sum over $\mu_i^+$ and use
$$\sum_{|\mu_i^+|=d} \frac{\chi_{\nu_i^+}(\mu_i^+)\chi_{\nu_{i+1}^-}(\mu_i^+)}{z_{\mu_i^+}}
= \delta_{\nu_i^+, \nu_{i+1}^-}.$$
We get
\begin{eqnarray*}
&& Z(\lambda) \\
& = & \prod_{i\in \bZ_k} \sum_{\nu_i^{+}}
e^{\sqrt{-1}\left(\kappa_{\nu_i^+} \frac{u^-_i}{u^+_i} + \kappa_{\nu_{i-1}^+}\frac{u_i^+}{u_i^-}\right)\lambda/2}
\cW_{\nu_i^+, \nu_{i-1}^+} ((-1)^{s_i}t_i)^{|\nu^+_i|} \\
& = & \prod_{i\in \bZ_k} \sum_{\nu_i}
e^{\sqrt{-1} \kappa_{\nu_i} \left(\frac{u^-_i}{u^+_i} + \frac{u_{i+1}^+}{u_{i+1}^-}\right)\lambda/2}
\cW_{\nu_i, \nu_{i-1}} ((-1)^{s_i}t_i)^{|\nu_i|} \\
& = & \prod_{i\in \bZ_k} \sum_{\nu_i}
e^{\sqrt{-1} \kappa_{\nu_i}s_i \lambda/2}
\cW_{\nu_i, \nu_{i-1}} ((-1)^{s_i}t_i)^{|\nu_i|}.
\end{eqnarray*}
In the last equality we have used (\ref{eqn:si}).
\end{proof}

\section{Examples}
\label{sec:Examples}

In Theorem \ref{thm:Main},
it suffices to know the homology classes (to determine $t_i$) and the self intersection numbers $s_i$.
In Iqbal's original version,
one has to know the details of the weights of the torus action at a vertex to write
the vertex operator $\cV$ in the form $T^kS^{-1}T^l$.
So our result is simpler than what Iqbal has conjectured,
as can be seen from the following examples.

\subsection{The resolved conifold case}

This case has been treated in \cite{Zho5}.
In this case,
we need the chemistry of two-colored labelled graphs.

\subsection{The local $\bP^1 \times \bP^1$ case}

Consider the following $T^2$-action on $\bP^1 \times \bP^1$:
$$(e^{\sqrt{-1}\vartheta_1}, e^{\sqrt{-1}\vartheta_2}) \cdot ([z_0:z_1], [w_0:w_1])
= [e^{\sqrt{-1}\vartheta_1}z_0:z_1], [e^{\sqrt{-1}\vartheta_2}w_0: w_1]).$$
Then $\bP^1\times \bP^1$ has four fixed points:
\begin{align*}
p_1 & = ([0:1], [0:1]), & p_3 & = ([0:1],[1:0]), \\
p_2 & = ([1:0],[0:1]), & p_4 & = ([1:0], [1:0]).
\end{align*}
There are four balloons $l_{p_{i}p_{i+1}}$, $i \in \bZ_4$.
It is not hard to see that
\begin{eqnarray*}
&& s_i = 0
\end{eqnarray*}
for all edges.
Furthermore,
in $H_2(\bP^1 \times \bP^1, \bZ)$,
we have
\begin{align*}
[l_{p_1p_2}] & = [l_{p_3p_4}], & [l_{p_2p_3}] & = [l_{p_4p_1}],
\end{align*}
hence
\begin{align*}
t_0 & = t_2, & t_2 & = t_3.
\end{align*}
The fixed points and the balloons can be put into the following picture:
$$
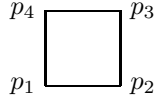
\xy
(10,0); (20,0), **@{-}; (7,0)*+{p_{4}}; (23,0)*+{p_{3}};
(10,-10); (20,-10), **@{-}; (7,-10)*+{p_{1}}; (23,-10)*+{p_{2}};
(10,0); (10,-10), **@{-};
(20,0); (20,-10), **@{-};
\endxy
$$
The prediction by Theorem \ref{thm:Main} is
\begin{eqnarray*}
&& Z(\lambda)
= \prod_{i\in \bZ_4} \sum_{\nu_i}
\cW_{\nu_i, \nu_{i-1}} t_i^{|\nu_i|}.
\end{eqnarray*}
This is exactly Iqbal's conjecture in this case (cf. \cite{Iqb}, (77)).

\subsection{The local $\bP^2$ case}
Consider the following $T^3$-action on $\bP^2$:
$$(e^{\sqrt{-1}\vartheta_1}, e^{\sqrt{-1}\vartheta_2}) \cdot [z_0:z_1:z_2]
= [z_0:e^{\sqrt{-1}\vartheta_1}z_1: e^{\sqrt{-1}\vartheta_2}z_2].$$
Then $\bP^2$ has three fixed points:
\begin{align*}
p_0 & = [1:0:0], & p_1 & = [0:1:0], & p_2 & = [0:0:1].
\end{align*}
Since the homology class of  each $l_{p_ip_{i+1}}$ is the same hyperplane class $H \in H_2(\bP^2, \bZ)$,
we have:
\begin{align*}
t_1 & = t_2 = t_3 = t, & s_1 = s_2 = s_3 = 1.
\end{align*}
The fixed points and balloons of $\bP^2$ can be put in the following picture:
$$
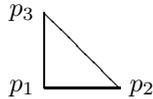
\xy
(10,0); (20,-10), **@{-};  (7, 0)*+{p_3};
(10,-10); (20,-10), **@{-}; (7,-10)*+{p_1}; (23,-10)*+{p_2};
(10,0); (10,-10), **@{-};
\endxy$$
Theorem \ref{thm:Main} predicts:
\begin{eqnarray*}
Z(\lambda)
& = & \prod_{i\in \bZ_3} \sum_{\nu_i}
e^{\sqrt{-1} \kappa_{\nu_i} \lambda/2}
\cW_{\nu_i, \nu_{i-1}} (-t)^{|\nu_i|} \\
& = & \sum_{\nu_1, \nu_2, \nu_3}\cW_{\nu_1, \nu_3}\cW_{\nu_2, \nu_1}\cW_{\nu_2, \nu_2}
q^{(\kappa_{\nu_1} + \kappa_{\nu_2} + \kappa_{\nu_3})/2}
(-t)^{|\nu_1| + |\nu_2| + |\nu_3|},
\end{eqnarray*}
where
$$q = e^{\sqrt{-1}\lambda}.$$
This is exactly Iqbal's conjecture in this case (cf. \cite{Iqb}, (36)).
This was obtained earlier by Aganagic, Marino, and Vafa \cite{Aga-Mar-Vaf} using different method.

\subsection{Local $\cB_1$ case}
Now we blow up $\bP^2$ in last subsection at $p_2$ and denote the resulting surface by $\cB_1$.
Then $p_2$ is replaced by an exception divisor $E$.
The $T$-action on $\bP^2$ lifts to $\cB_1$.
Now we have four fixed points and four balloons.
The original $p_1$ and $p_3$ are still fixed points.
We will rename $p_3$ as $p_4$.
There are two fixed points on the exceptional divisor $E$:
the intersection points of he strict transforms of the lines $l_{p_1p_2}$ and $l_{p_3p_2}$,
denoted by $p_2$ and $p_3$, respectively.
One can view $\cB_1$ as a $\bP^1$-fiber bundle over $E$,
then $H_2(\cB_1, \bZ)$ is generated by the class $B$ of $E$ and the class $F$ of a fiber.
It is well-known that
\begin{align*}
B^2 & = -1, & B \cdot F & = 1, & F^2 & = 0,
\end{align*}
and
\begin{align*}
[l_{p_1p_2}] & = [l_{p_3p_4}] = F, &  [l_{p_2p_3}] & = B, & [l_{p_4p_1}] & = B + F.
\end{align*}
It follows that
\begin{align*}
s_1 & = 0, & s_2 & = -1, & s_3 & = 0, & s_4 & = 1,
\end{align*}
and
\begin{align*}
t_3 & = t_1, & t_4 & = t_1t_2.
\end{align*}
The corresponding polygon is
$$
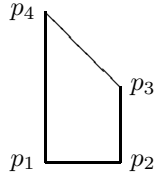
\xy
(10,0); (20,-10), **@{-};  (7, 0)*+{p_4};
(10,0); (10,-20), **@{-};
(10,-20); (20,-20), **@{-}; (7,-20)*+{p_1};
(20, -20); (20, -10), **@{-}; (23, -20)*+{p_2}; (23,-10)*+{p_3};
\endxy$$
Theorem \ref{thm:Main} predicts:
\begin{eqnarray*}
&& Z(\lambda)
= \sum_{\nu_1, \dots, \nu_4} \prod_{i \in \bZ_4}
\cW_{\nu_i, \nu_{i-1}} \cdot
q^{(\kappa_{\nu_4} - \kappa_{\nu_2})/2}
(-1)^{|\nu_4| - |\nu_2|} t_1^{|\nu_1| + |\nu_3| + |\nu_4|}
t_2^{|\nu_2|+|\nu_4|}.
\end{eqnarray*}
This is exactly Iqbal's conjecture for $\cB_1$ (\cite{Iqb}, (50)).

\subsection{Local $\cB_2$ case}
Now we consider the surface $\cB_2$ obtained from $\bP^2$ blown up at both $p_2$ and $p_3$.
The associated graph is given below:
$$
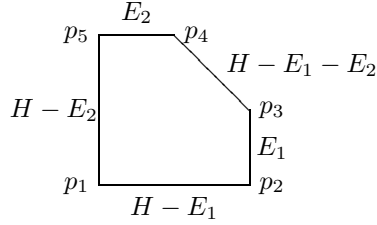
\xy
(10,0); (20,0), **@{-};  (7, 0)*+{p_5}; (23, 0)*+{p_4};
(10,0); (10,-20), **@{-};
(10,-20); (30,-20), **@{-}; (7,-20)*+{p_1};
(30, -20); (30, -10), **@{-}; (33, -20)*+{p_2}; (33,-10)*+{p_3};
(20, 0); (30, -10), **@{-};
(33, -15)*+{E_1}; (20, -23)*+{H-E_1}; (4, -10)*+{H-E_2}; (15, 3)*+{E_2}; (37, -4)*+{H-E_1-E_2};
\endxy$$
Here we have indicated the holomogy classes of the balloons,
where $H$ is the class of the strict transform of $l_{p_1p_2}$ plus the class of the exceptional divisor $E_1$.
Following Iqbal \cite{Iqb},
we will take $H, E_1, E_2$ as a basis of $H_2(\cB_2)$,
denote by $t_H$, $t_{E_1}$, $t_{E_2}$ the corresponding element in the Novikov ring.
It is clear that
\begin{align*}
t_1 & = t_Ht_{E_1}^{-1}, & t_2 & = t_{E_1}, & t_3 & = t_Ht_{E_1}^{-1}t_{E_2}^{-1}, &
t_4 & = t_{E_2}, & t_5 & = t_Ht_{E_2}^{-1}.
\end{align*}
From
\begin{align*}
E_1^2 & = E_2^2 = -1, & E_1 \cdot E_2 & = H\cdot E_1 = H \cdot E_2 = 0, & H^2 = 1,
\end{align*}
one gets
\begin{align*}
s_1 & = s_5 = 0, & s_2 & = s_3 =  s_4 = - 1.
\end{align*}
Therefore Theorem \ref{thm:Main} predicts:
\begin{eqnarray*}
 Z(\lambda)
& = & \sum_{\nu_1, \dots, \nu_5}
\prod_{i \in \bZ_5} \cW_{\nu_i, \nu_{i-1}} \cdot
q^{-(\kappa_{\nu_2} + \kappa_{\nu_3} + \kappa_{\nu_4})/2}
(-1)^{|\nu_2| + |\nu_3| + |\nu_4|} \\
&& \cdot t_H^{|\nu_1| + |\nu_3| + |\nu_5|} t_{E_1}^{-|\nu_1|+|\nu_2|-|\nu_3|}
t_{E_2}^{-|\nu_3|+|\nu_4|-|\nu_5|}.
\end{eqnarray*}
This is Iqbal's conjecture for $\cB_2$ (\cite{Iqb}, (64)).

\subsection{Local $\cB_3$ case}
Consider the surface $\cB_3$ obtained from $\bP^2$ by blowing up $p_1, p_2, p_3$.
The corresponding polygon is as follows:
$$
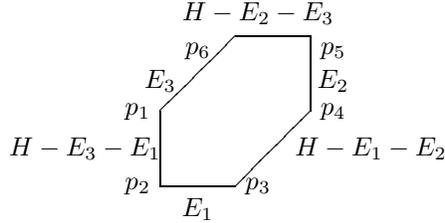
\xy
(10,-10); (20,0), **@{-};  (7, -10)*+{p_1}; (33, -2)*+{p_5}; (15, -2)*+{p_6};
(10,-10); (10,-20), **@{-};
(10,-20); (20,-20), **@{-}; (7,-20)*+{p_2};
(20, -20); (30, -10), **@{-}; (23, -20)*+{p_3}; (33,-10)*+{p_4};
(30, 0); (30, -10), **@{-}; (20, 0); (30, 0), **@{-};
(0, -15)*+{H-E_3-E_1}; (15, -23)*+{E_1}; (38, -15)*+{H-E_1-E_2};
(33, -6)*+{E_2}; (10, -6)*{E_3}; (23, 3)*+{H-E_2-E_3};
\endxy$$
Following Iqbal \cite{Iqb},
we take $H, E_1, E_2, E_3$ as a basis of $\cB_3$,
where $E_j$ is the exceptional divisor obtained by blowing up $p_j$,
$H$ is the class of the strict transform of $l_{p_1p_2}$ plus $E_1$ and $E_2$.
Denote by $t_H, t_{E_1}, t_{E_2}$, $t_{E_3}$ the corresponding elements in the Novikov ring of $\cB_3$.
Then we have ($j \in \bZ_3$):
\begin{align*}
t_{2j-1} & = t_Ht_{E_j}^{-1}t_{E_{j-1}}^{-1}, &
t_{2j} & = t_{E_j}.
\end{align*}
From
\begin{align*}
E_j^2 & = -1, & E_j \cdot E_{j+1} &= H \cdot E_j = 0, & H^2 = 1,
\end{align*}
we get
\begin{align*}
s_1 = s_2 = s_3 =  s_4 = s_5 = s_6 = - 1.
\end{align*}
The prediction by Theorem \ref{thm:Main} is
\begin{eqnarray*}
Z(\lambda)
& = & \sum_{\nu_1, \dots, \nu_6} \prod_{i \in \bZ_6}
\cW_{\nu_i, \nu_{i-1}}\cdot
q^{\sum_{i\in \bZ_6} \kappa_{\nu_i}/2}
(-1)^{\sum_{i \in \bZ_6} |\nu_i|} \\
&& \cdot  t_H^{|\nu_1| + |\nu_3| + |\nu_5|}
\prod_{j=1}^3 t_{E_j}^{-|\nu_{2-j}|+|\nu_{2j}|-|\nu_{2j+1}|}.
\end{eqnarray*}
This is Iqbal's conjecture for $\cB_3$ (\cite{Iqb}, (72)).

\end{document}